\documentclass[11pt,letterpaper,reqno]{amsart}
\usepackage{amsmath,amsthm,amsfonts,amssymb,mathtools,algpseudocode,algorithm}
\usepackage{comment,bookmark,bm,color,float,caption,subcaption}

\hypersetup{pdfstartview={FitH}}

\newtheorem{Theorem}{{\bf Theorem}}[section]

\newtheorem{Proposition}[Theorem]{{\bf Proposition}}

\numberwithin{equation}{section}

\newcommand{\C}{\mathbb{C}}
\newcommand{\R}{\mathbb{R}}

\newcommand{\calO}{\mbox{\Large $\mathcal{O}$}}
\newcommand{\calo}{\mathcal{O}}
\newcommand{\calB}{\mathcal{B}}
\newcommand{\calP}{\mathcal{P}}

\newcommand{\off}{\textup{off}}
\newcommand{\diag}{\text{diag}}

\newcommand{\doublewidetilde}[1]{\widetilde{\raisebox{0pt}[0.85\height]{$\widetilde{#1}$}}}
\newcommand{\doubletilde}[1]{\tilde{\raisebox{0pt}[0.85\height]{$\tilde{#1}$}}}

\begin{document}

\title[On the block Eberlein diagonalization method]{On the block Eberlein diagonalization method}
%\author{One author\thanks{some info} \and Another author\thanks{more info}}
\author{Erna Begovi\'{c}~Kova\v{c}}\thanks{\textsc{Erna Begovi\'{c} Kova\v{c},
University of Zagreb Faculty of Chemical Engineering and Technology, Maruli\'{c}ev trg 19, 10000 Zagreb, Croatia},
\texttt{ebegovic@fkit.unizg.hr}}
\author{Ana Perkovi\'{c}}\thanks{\textsc{Ana Perkovi\'{c},
University of Zagreb Faculty of Chemical Engineering and Technology, Maruli\'{c}ev trg 19, 10000 Zagreb, Croatia},
\texttt{aperkov@fkit.unizg.hr}}

\thanks{This work has been supported in part by Croatian Science Foundation under the project UIP-2019-04-5200.}
%\date{\today}

\renewcommand{\subjclassname}{\textup{2020} Mathematics Subject Classification}

\begin{abstract}
The Eberlein diagonalization method is an iterative Jacobi-type method for solving the eigenvalue problem of a general complex matrix. In this paper, we develop a block version of the Eberlein method. We prove the global convergence of our block algorithm and present several numerical examples.
\end{abstract}

\maketitle

\section{Introduction}

The Eberlein diagonalization method, introduced in ~\cite{Eberlein62}, is a Jacobi-type method for solving the eigenvalue problem of a general matrix. For a given matrix $A\in\C^{n\times n}$, the Eberlein algorithm is an iterative process of the form
\begin{equation}\label{Eber-elements}
A^{(k+1)}=T_k^{-1}A^{(k)}T_k, \quad A^{(0)}=A, \quad k\geq0.
\end{equation}
Each transformation $T_k$ in~\eqref{Eber-elements} is a product of a plane rotation $R_k$, chosen to annihilate the pivot element of the Hermitian part of $A^{(k)}$, $B^{(k)}=\left( A^{(k)}+\left(A^{(k)}\right)^*\right)/2$, and a norm reducing transformation $S_k$;
$$T_k=R_kS_k.$$
Compared to the other Jacobi-type methods, the importance of the Eberlein method lies in the fact that it can be applied to any matrix $A$. That is, the starting matrix $A$ does not need to have any specific structure or properties.

Convergence of this method was studied in ~\cite{Veselic76,Hari82,PupHari98,BePe24}. For different pivot strategies, it was shown that the iterations~\eqref{Eber-elements} converge in the sense that the sequence $(B^{(k)},k\geq0)$ converges to a diagonal matrix, while $(A^{(k)},k\geq0)$ converges to a normal matrix $\Lambda$. If all the eigenvalues of $A$ have different real parts, then $\Lambda$ is a diagonal matrix. Otherwise, $\Lambda$ is permutation-similar to a block diagonal matrix, such that the sizes of its blocks correspond to the multiplicities of the real parts of the eigenvalues of $A$.

In this paper, we present the block Eberlein method. To the best of our knowledge, this is the first block variant of the Eberlein algorithm. This makes an important step in the development of the Jacobi-type methods. The use of the matrix blocks, instead of the elements, can improve the efficiency of the algorithm, since a block algorithm can exploit the cache hierarchy of modern computers. 

We formulate the iterative procedure on block matrices and determine the choice of the transformation matrices. Then, we prove the convergence of our block method. The convergence results are in line with those for the element-wise method. Regarding the pivot orderings, we work with the generalized serial pivot orderings from~\cite{BeHa17}. This is a very wide class of pivot orderings that includes the most common serial orderings (row-wise and column-wise), as well as many other orderings. 

The paper is structured as follows. In Section~\ref{sec:prelim} we give preliminary results and introduce the notation. We describe our algorithm in Section~\ref{sec:method} and prove its convergence in Section~\ref{sec:convergence}. Then, in Section~\ref{sec:numerical}, we present different numerical examples. We end the paper with a short conclusion in Section~\ref{sec:conclusion}.

\section{Preliminaries}\label{sec:prelim}

We denote the block matrices by boldface capital letters, e.g., $\mathbf{A}$, $\mathbf{B}$. The block partition of an $n\times n$ matrix $\mathbf{A}$ is determined by an integer partition
\begin{equation}\label{partition_pi}
    \pi=(n_1, n_2,\ldots, n_m),
\end{equation}
where $n_i\geq 1$, for all $1\leq i\leq m$, and $n_1+n_2+\cdots+n_m=n$. Then,
$$\mathbf{A}=\left[\begin{array}{cccc}
      A_{11} &A_{12} & \ldots & A_{1m} \\
      A_{21} & A_{22} & \ldots & A_{2m} \\
      \vdots & \vdots & \ddots & \vdots \\
      A_{m1} & A_{m2} &\ldots &A_{mm}
    \end{array} \right]
    \begin{array}{r}
        n_1\\ n_2\\ \vdots \\ n_m
    \end{array}.$$
The dimension of the matrix block $A_{ij}$ is $n_i\times n_j$, for all $1\leq i,j\leq m$. Obviously, if 
\begin{equation}\label{pi1}
\pi=\pi_1\coloneqq(1,1,...1),
\end{equation}
then the block $A_{ij}$ is actually one element, $a_{ij}$. 

An \emph{elementary block matrix} $\mathbf{E}_{pq}$ with block partition~\eqref{partition_pi} is a block matrix that differs from the identity in only four blocks, those at the intersection of the $p$th and $q$th block row and column. We have
\begin{equation}\label{elemblockm}
    \mathbf{E}_{pq}=\left[\begin{array}{ccccc}
      I &  &  & & \\
       & E_{pp} &  & E_{pq} & \\
       & & I & & \\
       & E_{pq} &  & E_{qq} & \\
       &  &  & &  I\\
    \end{array} \right]
    \begin{array}{r}
        \\ n_p\\ \\ n_q \\  \\
    \end{array}.
\end{equation}
For
\begin{equation*}
    \widehat{\mathbf{E}}_{pq}=\left[\begin{array}{cc}
        E_{pp} & E_{pq} \\
        E_{qp} & E_{qq}
    \end{array}\right],
\end{equation*}
the function that maps $\widehat{\mathbf{E}}_{pq}$ to an elementary block matrix $\mathbf{E}_{pq}$ with partition $\pi$ is denoted by $\mathcal{E}_{\pi}$. We write
$$\mathbf{E}_{pq}=\mathcal{E}_{\pi}(p,q,\widehat{\mathbf{E}}_{pq}).$$

\emph{Off-norm} of a matrix $A$ is defined as the Frobenius norm of its off-diagonal part,
$$\off(A)=\|A-\diag(A)\|_F.$$
Then, off-norm of 
a block matrix $\mathbf{A}$ with the partition $\pi$ as in~\eqref{partition_pi} can be written as
$$\off^2(\mathbf{A})=\sum_{\substack{i,j=1 \\ i\neq j}}^m\|A_{ij}\|_F^2+\sum_{i=1}^m\off^2(A_{ii}).$$
If $\off(\mathbf{A})=0$, then $\mathbf{A}$ is a diagonal matrix.

Next, we recall the element-wise Eberlein diagonalization method. In the element-wise iterations~\eqref{Eber-elements}, both $R_k$ and $S_k$ are elementary matrices, differing from the identity in only one principal $2\times2$ submatrix,
\begin{equation}\label{eq:RS}
    \widehat{R}_k=\begin{bmatrix}
        \cos \varphi_{k} & -e^{\imath\alpha_k} \sin \varphi_k\\
        e^{-\imath\alpha_k}\sin \varphi_k & \cos \varphi_k
    \end{bmatrix}, \quad
    \widehat{S}_k=\begin{bmatrix}
        \cosh\psi_k &  -\imath e^{\imath\beta_k}\sinh \psi_k\\
         \imath e^{-\imath\beta_k}\sinh \psi_k & \cosh \psi_k
    \end{bmatrix},
\end{equation}
respectively, where $\imath$ stands for the imaginary unit.
The position of the submatrices $\widehat{R}_k$ and $\widehat{S}_k$ within $R_k$ and $S_k$ is determined by the \emph{pivot pair} $(p,q)=(p(k),q(k))$. We have $R_k=\mathcal{E}_{\pi_1}(p,q,\widehat{R}_k)$ and $S_k=\mathcal{E}_{\pi_1}(p,q,\widehat{S}_k)$, for $\pi_1$ like in~\eqref{pi1}.
When the behavior of the Eberlein method is examined, apart from observing the sequence of matrices $(A^{(k)},k\geq0)$ generated by~\eqref{Eber-elements}, two other matrix sequences are important to look at. One is a sequence of the already mentioned Hermitian parts of $A^{(k)}$, denoted by $(B^{(k)},k\geq0)$. The other one is the sequence obtained by applying the operator $C$, 
\begin{equation}\label{eq:C}
C\big(A^{(k)}\big)\coloneqq A^{(k)}\big(A^{(k)}\big)^*-\big(A^{(k)}\big)^*A^{(k)},
\end{equation}
which, in a way, measures the distance of $A^{(k)}$ from the set of normal matrices. 

In the block Eberlein method, block transformations will be elementary block matrices of the form~\eqref{elemblockm}, related to those in~\eqref{eq:RS}. We are going to observe the sequence of block matrices $(\mathbf{A}_k, k\geq0)$ generated by our block Eberlein algorithm and the sequence of their Hermitian parts $(\mathbf{B}_k, k\geq0)$. The operator $C$ will be applied on block matrices in the same way as in the element-wise case.

\section{Block Eberlein diagonalization algorithm}\label{sec:method}

Now we are going to describe the block Eberlein algorithm.
Let $\mathbf{A}\in\C^{n\times n}$ be an arbitrary block matrix with partition $\pi$ as in~\eqref{partition_pi}. The block Eberlein method is the iterative process 
\begin{equation}\label{blockEber_process}
\mathbf{A}^{(k+1)}=\mathbf{T}_k^{-1}\mathbf{A}^{(k)}\mathbf{T}_k, \quad k\geq0,
\end{equation}
where $\mathbf{A}^{(0)}=\mathbf{A}$, and
\begin{equation}\label{transformationT}
\mathbf{T}_k=\mathbf{R}_k\mathbf{S}_k, \quad k\geq 0,
\end{equation}
are non-singular elementary block matrices. The partition $\pi$ is the same for all matrices from the relations~\eqref{blockEber_process} and~\eqref{transformationT}. It is assumed to be fixed throughout the process, so we omit it from the notation. However, it would be possible to have an adaptive partition that is changing throughout the iterations. 
Transformations $\mathbf{T}_k=\mathcal{E}_{\pi}(p(k),q(k),\widehat{\mathbf{T}}_{p(k)q(k)})$ are elementary block transformations. 
Precisely, $\mathbf{R}_k=\mathcal{E}_{\pi}(p(k),q(k),\widehat{\mathbf{R}}_{p(k)q(k)})$ are unitary elementary block transformations chosen to diagonalize the  pivot submatrix of the Hermitian part of $\mathbf{A}^{(k)}$,
\begin{equation}\label{eq:blockB}
    \mathbf{B}^{(k)}=\frac{1}{2}\left(\mathbf{A}^{(k)}+\big(\mathbf{A}^{(k)}\big)^*\right),
\end{equation}
while $\mathbf{S}_k=\mathcal{E}_{\pi}(p(k),q(k),\widehat{\mathbf{S}}_{p(k)q(k)})$ are nonsingular nonunitary elementary block transformations that reduce the Frobenius norm of $\mathbf{A}^{(k)}$. 
The index pair $(p,q)=(p(k),q(k))$ determines the $k$th \emph{pivot block}. For the sake of simplicity of notation, when $k$ is implied, we will omit it and write $(p,q)$.

The process~\eqref{blockEber_process} can be written with an intermediate step, $$\mathbf{A}^{(k)}\to  \widetilde{\mathbf{A}}^{(k)} \to \mathbf{A}^{(k+1)},$$ where
\begin{align}
    \widetilde{\mathbf{A}}^{(k)}&=\mathbf{R}_k^* \mathbf{A}^{(k)} \mathbf{R}_k, \label{eq:blockpart1} \\ 
    \mathbf{A}^{(k+1)}&=\mathbf{S}_k^{-1} \widetilde{\mathbf{A}}^{(k)} \mathbf{S}_k,\quad k\geq 0.\label{eq:blockpart2}
\end{align}
We provide the details on the steps~\eqref{eq:blockpart1} and~\eqref{eq:blockpart2} in the Subsections~\ref{sec:Rk} and~\ref{sec:Sk}, respectively.
The whole procedure for the block Eberlein method is given in Algorithm~\ref{agm:blockeberlein}.
Clearly, if $\pi=\pi_1$ from the relation~\eqref{pi1}, then the block Eberlein method comes down to the element-wise Eberlein method. 

\begin{algorithm}
\caption{Block Eberlein method}
\label{agm:blockeberlein}
\begin{algorithmic}
\State \textbf{Input:} block matrix $\mathbf{A}$
\State \textbf{Output:} block matrix $\Lambda=\mathbf{A}^{(K)}$, block elementary matrix $\mathbf{T}_K$
\State $\mathbf{A}^{(0)}=\mathbf{A}$, $\mathbf{T}_0=I$
\State $k=0$
\Repeat
\State Choose block pivot pair $(p,q)$ according to the pivot strategy.
\State Find $\widehat{\mathbf{R}}^{(k)}_{pq}$ that diagonalizes the Hermitian matrix $\widehat{\mathbf{B}}_{pq}^{(k)}$. \Comment{complex Jacobi algorithm}
\State $\mathbf{R}_k=\mathcal{E}(p,q,\widehat{\mathbf{R}}^{(k)}_{pq})$
\State $\widetilde{\mathbf{A}}^{(k)}=\mathbf{R}_k^*{\mathbf{A}}^{(k)}\mathbf{R}_k$
\State Find $\widehat{\mathbf{S}}^{(k)}_{pq}$ which reduces the Frobenius norm of $\widetilde{\mathbf{A}}^{(k)}$. \Comment{Algorithm~\ref{agm:Sblock}}
\State $\mathbf{S}_k=\mathcal{E}(p,q,\widehat{\mathbf{S}}^{(k)}_{pq})$
\State ${\mathbf{A}}^{(k+1)}=\mathbf{S}_k^{-1}\widetilde{\mathbf{A}}^{(k)}\mathbf{S}_k$
\State $\mathbf{T}_{k+1}=\mathbf{T}_k \mathbf{R}_k \mathbf{S}_k$
\State $k=k+1$
\Until{convergence}
\State $\mathbf{A}^{(K)}=\mathbf{A}^{(k)}$, $\mathbf{T}_K=\mathbf{T}_k$
\end{algorithmic}
\end{algorithm}

\subsection{Unitary block transformations}\label{sec:Rk}

In the $k$th step of the method, the unitary elementary block transformation $\mathbf{R_k}$ is chosen to diagonalize the pivot submatrix of $\mathbf{B}^{(k)}$ from~\eqref{eq:blockB}. We have
$$\left(\widehat{\mathbf{R}}^{(k)}_{pq}\right)^* \widehat{\mathbf{B}}^{(k)}_{pq} \widehat{\mathbf{R}}^{(k)}_{pq} = \left[\begin{array}{cc}
    R^{(k)}_{pp} & R^{(k)}_{pq} \\
    R^{(k)}_{qp} & R^{(k)}_{qq}
\end{array}\right]^*\left[\begin{array}{cc}
    B^{(k)}_{pp} & B^{(k)}_{pq} \\
    B^{(k)}_{qp} & B^{(k)}_{qq}
\end{array}\right]\left[\begin{array}{cc}
    R^{(k)}_{pp} & R^{(k)}_{pq} \\
    R^{(k)}_{qp} & R^{(k)}_{qq}
\end{array}\right]=\left[\begin{array}{cc}
    \Lambda^{(k+1)}_{pp} & 0 \\
    0 & \Lambda^{(k+1)}_{qq}
\end{array}\right],
$$
where $\widehat{\mathbf{R}}^{(k)}_{pq}$ and $\widehat{\mathbf{B}}^{(k)}_{pq}$ are $(n_p+n_q)\times(n_p+n_q)$ matrices
and $\Lambda^{(k+1)}_{pp}$ and $\Lambda^{(k+1)}_{qq}$ are diagonal matrices.
In order to determine $\widehat{\mathbf{R}}_{pq}^{(k)}$, we observe $\widehat{\mathbf{B}}_{pq}^{(k)}$ element-wise and apply the complex Jacobi method from~\cite{BeHa21}. Also, one can take $\widehat{\mathbf{B}}_{pq}^{(k)}$ as a $2\times2$ block matrix and apply the complex block Jacobi method from~\cite{BeHa24}. Instead of the Jacobi method, other diagonalization methods could be considered as well (see, e.g.,~\cite{Drmac21}).

Recall that for the standard element-wise Jacobi diagonalization method,
the sufficient convergence condition is the existence of a strictly positive uniform lower bound for the cosine of the rotation angle,~\cite{FoHe60}. A generalization of that condition to the block Jacobi method is the existence of such bound for the singular values of the diagonal blocks of transformation matrices,~\cite{Drmac09}. 

Unitary block matrix $\mathbf{U}$ with block partition $\pi=(n_1,n_2)$, $n_1,n_2\geq1$, $n=n_1+n_2$,
$$U=\begin{bmatrix}
U_{11} & U_{12} \\
U_{21} & U_{22} \\
\end{bmatrix},$$
such that the singular values of its diagonal blocks can be bounded from below by a function of dimension is called \emph{UBC (uniformly bounded cosine) matrix}, defined in~\cite{Drmac09}. 
In the same paper it was proven that, for any $n\times n$ unitary block matrix $\mathbf{U}$ with block partition $\pi=(n_1,n_2)$, there is a permutation $P$ such that, for $\mathbf{U}'=\mathbf{U}P$, we have
$\sigma_{\min}(\mathbf{U}'_{11})\geq\gamma_{\pi}>0,$
where $\gamma_{\pi}$ is a constant depending only on $\pi$. Moreover, it was shown in~\cite{Hari15} that there is a lower bound
$$\sigma_{\min}(\mathbf{U}'_{11})\geq\gamma_n>0,$$
depending only on $n$.
Therefore, every unitary elementary block matrix can be transformed into a UBC matrix using an appropriate permutation $P$. In conclusion, we can take block unitary transformations $\mathbf{R}_k$ in~\eqref{transformationT} to be UBC transformations. This will be important for the convergence results in Section~\ref{sec:convergence}.

\subsection{Norm-reducing block transformations}\label{sec:Sk}

The goal of the elementary block transformation $\mathbf{S}_k$ is to reduce the Frobenius norm of the matrix $\mathbf{A}^{(k)}$, that is, $\widetilde{\mathbf{A}}^{(k)}$ obtained in~\eqref{eq:blockpart1}. This is achieved by reducing the Frobenius norm of the pivot block columns and rows. Computing
$\widehat{\mathbf{S}}^{(k)}_{pq}$ is more demanding than computing $\widehat{\mathbf{R}}^{(k)}_{pq}$. In order to obtain $\widehat{\mathbf{R}}^{(k)}_{pq}$, it is enough to consider the four blocks, that is, the pivot submatrix $\widehat{\mathbf{B}}^{(k)}_{pq}$ of $\mathbf{B}^{(k)}$. On the other hand, in accordance with the element-wise Eberlein method, obtaining $\widehat{\mathbf{S}}^{(k)}_{pq}$ requires $4m-4$ blocks, that is, $p$th and $qth$ block column and row.

For a fixed pivot pair $(p,q)$, we construct the core algorithm for finding ${\mathbf{S}}_k$ based on the norm-reducing transformations in the element-wise Eberlein algorithm, although, the reduction of the Frobenius norm can be achieved in more than one way. We observe the iterative process
\begin{equation}\label{eq:coreprocess}
\doublewidetilde{\mathbf{A}}^{(l+1)}=\widetilde{\mathbf{S}}_{l}^{-1}\doublewidetilde{\mathbf{A}}^{(l)}\widetilde{\mathbf{S}}_{l}, \quad l\geq 0,
\end{equation}
where $\doublewidetilde{\mathbf{A}}^{(0)}=\widetilde{\mathbf{A}}^{(k)}$, and $\widetilde{\mathbf{S}}_{l}\in\C^{n\times n}$ are of the form
\begin{equation}\label{Sl}
\widetilde{\mathbf{S}}_l=\left[
    \begin{array}{ccccccccccc}
      1 &  &  &  &  &  &  &  &  &  &  \\
       & \ddots &  &  &  &  &  &  &  &  &  \\
       &  & 1 &  &  &  &  &  &  &  &  \\
       &  &  & \cosh\psi_l &  &  &  & -\imath e^{\imath\beta_l}\sinh\psi_l &  &  &  \\
       &  &  &  & 1 &  &  &  &  &  &  \\
       &  &  &  &  & \ddots &  &  &  &  &  \\
       &  &  &  &  &  & 1 &  &  &  &  \\
       &  &  & \imath e^{-\imath\beta_l}\sinh\psi_l &  &  &  & \cosh\psi_l &  &  &  \\
       &  &  &  &  &  &  &  & 1 &  &  \\
       &  &  &  &  &  &  &  &  & \ddots &  \\
       &  &  &  &  &  &  &  &  &  & 1 \\
    \end{array}
  \right]
  \begin{array}{l}
     \\
     \\
     \\
     r \\
     \\
     \\
     \\
     s \\
     \\
     \\
     \\
     \end{array}.
\end{equation}
For $l\geq0$, the transformation angles are calculated from the relations
\begin{equation}\label{beta}
\tan\beta_l = -\frac{\text{Re}\big(\doubletilde{c}_{rs}^{(l)}\big)}{\text{Im}\big(\doubletilde{c}_{rs}^{(l)}\big)},
\end{equation}
where
$$C\Big(\doublewidetilde{\mathbf{A}}^{(l)}\Big)=\big(\doubletilde{c}_{rs}^{(l)}\big),$$
and
\begin{equation}\label{psi}
\tanh\psi_l =\frac{\text{Im}\big(t_{rs}^{(l)}d_{rs}^{(l)*}\big)-w_{rs}^{(l)}/2} {v_{rs}^{(l)}+2\big(|t_{rs}^{(l)}|^2+|d_{rs}^{(l)}|^2\big)},
\end{equation}
where
\begin{align*}
d_{rs}^{(l)} & = \doubletilde{a}_{rr}^{(l)}-\doubletilde{a}_{ss}^{(l)}, \\
t_{rs}^{(l)} & = \big(\doubletilde{a}_{rs}^{(l)}+\doubletilde{a}_{rs}^{(l)}\big)\cos\beta_l -\imath\big(\doubletilde{a}_{rs}^{(l)}-\doubletilde{a}_{rs}^{(l)}\big)\sin\beta_l, \\
v_{rs}^{(l)} & = \sum_{\substack{i=1 \\ i\neq r,s}}^n |\doubletilde{a}_{ir}^{(l)}|^2+|\doubletilde{a}_{ri}^{(l)}|^2+|\doubletilde{a}_{is}^{(l)}|^2+|\doubletilde{a}_{si}^{(l)}|^2, \\
w_{rs}^{(l)} & = -\text{Re}\big(\xi_{rs}^{(l)}\big)\sin\beta_l +\text{Im}\big(\xi_{rs}^{(l)}\big)\cos\beta_l, \\
\xi_{rs}^{(l)} & = 2\sum_{\substack{i=1 \\ i\neq r,s}}^n \big(\doubletilde{a}_{ri}^{(l)}\doubletilde{a}_{si}^{(l)*}-\doubletilde{a}_{ir}^{(l)*}\doubletilde{a}_{is}^{(l)}\big).
\end{align*}
In the upper relations, $a^*$ stands for the complex conjugate of $a$. Such choice of transformations $\widetilde{\mathbf{S}}_l$ corresponds to the norm-reducing transformations for the element-wise case. It approximates the maximal norm reduction for each $(r,s)$. In particular, it follows from~\cite{Eberlein62} that
\begin{equation}\label{eq:delta}
\Delta_l:=\|\doublewidetilde{\mathbf{A}}^{(l)}\|_F^2 -\|\doublewidetilde{\mathbf{A}}^{(l+1)}\|_F^2  \geq\frac{1}{3}\frac{|\doubletilde{c}_{rs}^{(l)}|^2}{\|\mathbf{A}\|_F^2}, \quad l\geq0.
\end{equation}

Index pairs $(r,s)=(r(l),s(l))$, $r<s$, in~\eqref{Sl} are taken from the upper triangle of the pivot submatrix of $\widetilde{\mathbf{A}}^{(k)}$, which is determined by the pivot pair $(p,q)$.
To be precise, we have 
\begin{equation}\label{eq:index_rs}
\begin{aligned}
n_1+\cdots+n_{p-1}+1 & \leq r <s\leq n_1+\cdots+n_{p-1}+n_{p}, \\
& \textup{or} \\
n_1+\cdots+n_{p-1}+1 & \leq r \leq n_1+\cdots+n_{p-1}+n_{p}, \\
n_1+\cdots+n_{q-1}+1 & \leq s \leq  n_1+\cdots+n_{q-1}+n_{q},  \\
& \textup{or} \\
n_1+\cdots+n_{q-1}+1 & \leq r <s\leq n_1+\cdots+n_{q-1}+n_{q}. 
\end{aligned}
\end{equation}
Thus, the iterations~\eqref{eq:coreprocess} affect only the two pivot block columns and rows and transformations $\widetilde{\mathbf{S}}_l$ are of the form $\mathcal{E}_{\pi}(p,q,\widehat{\widetilde{\mathbf{S}}}_l)$.
Then, the submatrix $\widehat{\mathbf{S}}_{pq}$ is computed as 
$$\widehat{\mathbf{S}}_{pq}=\prod_l\widehat{\widetilde{\mathbf{S}}}_l.$$

Concerning the index $l$ from the upper product, in our implementation we take every possible position from~\eqref{eq:index_rs} exactly once. That is, we have 
\begin{equation}\label{indexl}
0\leq l<(n_p+n_q)(n_p+n_q-1)/2=\colon L.
\end{equation}
Then, $\mathbf{A}^{(k+1)}=\doublewidetilde{\mathbf{A}}^{(L-1)}$.
For the convergence properties it will only be important that 
$\mathbf{S}_k=\mathcal{E}_{\pi}(p,q,\widehat{\mathbf{S}}_{pq})$ reduces the Frobenius norm. If we take more sweeps over $l$, the reduction will be bigger, but the procedure will be more computationally exhausting. If we take only one pair from~\eqref{eq:index_rs}, the computation would be fast, but the norm reduction would be far from optimal, which would slower the convergence.

Algorithm~\ref{agm:Sblock} summarizes the discussion from this subsection.

\begin{algorithm}
\caption{Finding $\widehat{\mathbf{S}}^{(k)}_{pq}$}
\label{agm:Sblock}
\begin{algorithmic}
\State \textbf{Input:} block matrix $\widetilde{\mathbf{A}}^{(k)}$, pivot pair $(p,q)$
\State \textbf{Output:} submatrix $\widehat{\mathbf{S}}^{(k)}_{pq}$
\State $\doublewidetilde{\mathbf{A}}^{(0)}=\widetilde{\mathbf{A}}^{(k)}$, $\widehat{\mathbf{S}}^{(0)}_{pq}=I_{n_p+n_q}$
\State $l=0$
\Repeat
\State Choose pair $(r,s)$ from~\eqref{eq:index_rs}.
\State Find ${(n_p+n_q)\times (n_p+n_q)}$ block matrix $\widehat{\mathbf{S}}^{(l)}_{rs}$ using $\doublewidetilde{\mathbf{A}}^{(l)}$ and the relations~\eqref{beta} and~\eqref{psi}. 
\State $\mathbf{S}_{l}=\mathcal{E}(p,q,\widehat{\mathbf{S}}^{(l)}_{rs})$
\State $\doublewidetilde{\mathbf{A}}^{(l+1)}=\mathbf{S}_{l}^{-1}\doublewidetilde{\mathbf{A}}^{(l)}\mathbf{S}_{l}$
\State $\widehat{\mathbf{S}}^{(l+1)}_{pq}=\widehat{\mathbf{S}}^{(l)}_{pq}\widehat{\mathbf{S}}^{(l)}_{rs}$
\State $l=l+1$  
\Until{stopping criterion is satisfied}
\State $\widehat{\mathbf{S}}^{(k)}_{pq}= \widehat{\mathbf{S}}^{(l)}_{pq} $
\end{algorithmic}
\end{algorithm}

\subsection{Pivot orderings}

Just as any Jacobi-type algorithm, the block Eberlein algorithm depends on a \emph{block pivot ordering}, that is, the order of pivot blocks. In a block matrix $\mathbf{A}$ with block partition $\pi=(n_1,n_2,\ldots,n_m)$, possible pivot blocks are those from the upper triangle. We denote them by $\calP_m=\{(p,q) \ | \ 1\leq p<q\leq m\}$. 
A cyclic pivot ordering is a periodic ordering that repeatedly, in some prescribed order, takes all $M=\frac{m(m-1)}{2}$ pivot pairs from $\calP_m$. We denote the set of all cyclic ordering by $\calO(\calP_m)$ and say that $\calo\in\calO(\calP_m)$ is a cyclic pivot ordering.
It can be depicted by an $m\times m$ strictly upper triangular matrix $O$, such that
$$O=(o_{ij})=k \text{ if } (p(k),q(k))=(i,j), \quad \text{for } 1\leq i<j\leq m \text{ and }0\leq k<M.$$
The most common pivot orderings are the serial ones, row-wise and column-wise, where the pivot positions are taken cyclically row-by-row or column-by-column. 
Then, 
$$O_{\text{row}}=
\begin{bmatrix}
    \centering
    \begin{tabular}{cccc}
        * & 0 & 1 & 2 \\
        * & * & 3 & 4 \\
        * & * & * & 5 \\
        * & * & * & * \\
    \end{tabular}
\end{bmatrix} \quad \text{and} \quad O_{\text{col}}=
\begin{bmatrix}
    \centering
    \begin{tabular}{cccc}
        * & 0 & 1 & 3 \\
        * & * & 2 & 4 \\
        * & * & * & 5 \\
        * & * & * & * \\
    \end{tabular}
\end{bmatrix},
$$
respectively, represent row-wise ($\calo_{\text{row}}^{(4)}$) and column-wise ($\calo_{\text{col}}^{(4)}$) pivot ordering on a $4\times4$ block matrix.

We work with the \emph{generalized serial pivot orderings}. This class of block orderings, introduced in~\cite{BeHa17} and denoted by $\calB_{sg}$, is interesting because of its size. It includes the previously mentioned serial orderings, but also many others. To illustrate the size of this class of orderings, we can say that on the $4\times4$ matrices it covers $13/15\approx0.87$ of all possible pivot orderings~\cite{BeHa4x4parallel,BeHa4x4all}.

Let us define the set of block orderings $\calB_{sg}$. (One can check~\cite{BeHa17} for details.) First, note that an \emph{admissible transposition} of two adjacent pivot pairs $(i_k,j_k)$ and $(i_{k+1},j_{k+1})$ of a pivot ordering $\calo$ can be done if the indices $\{i_k,j_k,i_{k+1},j_{k+1}\}$ are all different. In the context of parallel strategies, transposition of the pivot pairs $(i_k,j_k)$ and $(i_{k+1},j_{k+1})$ is admissible if they can be used in parallel, at the same time.
If $\calo=(i_0,j_0),(i_1,j_1),\ldots,(i_{M-1},j_{M-1})$, then two pivot orderings $\calo,\calo'\in\calO(\calP_m)$ are
\begin{itemize}
\item[(i)] \emph{equivalent} $(\calo\sim\calo')$ if one can be obtained from the other by a finite set of admissible transpositions;
\item[(ii)] \emph{shift-equivalent} $(\calo \stackrel{s}{\sim} \calo')$ if $\calo'=(i_t,j_t),\ldots,(i_{M-1},j_{M-1}),(i_0,j_0),\ldots,(i_{t-1},j_{t-1})$, for some $1\leq t\leq m$;
\item[(iii)] \emph{weak-equivalent} $(\calo\stackrel{w}{\sim}\calo')$ if one can be obtained from the other by a finite set of equivalences or shift-equivalences
\item[(iv)] \emph{permutation-equivalent} $(\calo\stackrel{p}{\sim}\calo')$ if $\calo'=(\textup{q}(i_0),\textup{q}(j_0)),\ldots,(\textup{q}(i_{M-1}),\textup{q}(j_{M-1}))$, for some permutation $\textup{q}$ of length $M$.
\item[(v)] \emph{reverse} $(\calo'=\calo^{\gets})$ if $\calo'=(i_{M-1},j_{M-1}),\ldots,(i_1,j_1),(i_0,j_0)$.
\end{itemize}

Set
\begin{align*}
    \calB_c^{(m)}=\Big\{ \calo \in \calO(\calP_m) \mid \calo= &(1,2),(\tau_3(1),3),(\tau_3(2),3),\ldots,(\tau_m(1),m),\ldots\\
    &\ldots,(\tau_m(m-1),m), \quad \tau_j\in\Pi^{(1,j-1)}, 3\leq j\leq m \Big\}
\end{align*}
and
\begin{align*}
    \calB_r^{(m)}=\Big\{ \calo \in \calO(\calP_m) \mid \calo= &(m-1,m),(m-2,\tau_{m-2}(m-1)),(m-2,\tau_{m-2}(m)),\ldots\\
    &\ldots,(1,\tau_1(2)),\ldots, (1,\tau_1(m)), \quad \tau_i\in\Pi^{(i+1,m)}, 1\leq i\leq m-2 \Big\},
\end{align*}
where $\Pi^{(l_1,l_2)}$ stands for the set of all permutations of the set $\{l_1,l_1+1,l_1+2,\ldots,l_2\}$.
Pivot orderings from $\calB_c^{(m)}$ are derived from the column-wise ordering $\calo_{\text{col}}$. They take pivot blocks column-by-column, from left to right, but inside each column, pivot positions can be taken in an arbitrary order. Similarly, orderings from $\calB_r^{(m)}$ go row-by-row, from bottom to top, arbitrary inside each row. Then
$$\calB_{sp}^{(m)}=\calB_c^{(m)}\cup \overleftarrow{\calB}_c^{(m)} \cup \calB_r^{(m)}\cup \overleftarrow{\calB}_r^{(m)}$$
is the set of serial block pivot ordering with permutations. Our aimed set of the generalized block pivot orderings is an expansion of $\calB_{sp}^{(m)}$,
$$\calB^{(m)}_{sg}=\Big\{\calo\in\calO(\calP_m) \mid \calo\stackrel{w}{\sim}\calo'\stackrel{p}{\sim}\calo'' \text{ or } \calo\stackrel{p}{\sim}\calo'\stackrel{w}{\sim}\calo'',  \calo''\in\calB_{sp}^{(m)}\Big\}.$$
The orderings from $\calB^{(m)}_{sp}$ and $\calB^{(m)}_{sg}$ are very different from $\calo_{\text{row}}^{(m)}$ and $\calo_{\text{col}}^{(m)}$. For example
$$O_1=
\begin{bmatrix}
    \centering
    \begin{tabular}{ccccc}
        * & 0 & 2 & 4 & 7 \\
        * & * & 1 & 5 & 9 \\
        * & * & * & 3 & 6 \\
        * & * & * & * & 8 \\
        * & * & * & * & \\
    \end{tabular}
\end{bmatrix} \quad \text {and} \quad O_2=
\begin{bmatrix}
    \centering
    \begin{tabular}{ccccc}
        * & 7 & 2 & 0 & 6 \\
        * & * & 5 & 3 & 9 \\
        * & * & * & 8 & 4 \\
        * & * & * & * & 1 \\
        * & * & * & * & \\
    \end{tabular}
\end{bmatrix}$$
are matrix representations of the orderings $\calo_1\in\calB^{(5)}_{sp}$ and $\calo_2\in\calB^{(5)}_{sg}$, respectively.

\section{Convergence of the block Eberlein method}\label{sec:convergence}

After we have derived the block Eberlein method in the previous section, in this section we are going to prove that it is indeed convergent.  To that end, we will use several auxiliary results. First, in Theorem~\ref{tm:blockJacobi} we observe a different Jacobi-type process.

\begin{Theorem}\label{tm:blockJacobi}
Let $\mathbf{H}\in\C^{n\times n}$ be a Hermitian block matrix with partition $\pi=(n_1,n_2,\ldots,n_m)$, and let
$$\mathbf{H}^{(k+1)}=\mathbf{U}_k^*\mathbf{H}^{(k)}\mathbf{U}_k + \mathbf{M}^{(k)}, \quad \mathbf{H}^{(0)}=\mathbf{H}, \quad k\geq0,$$
with
$$\lim_{k\rightarrow\infty} \off\left(\mathbf{M}^{(k)}\right)=0.$$
If the pivot strategy is defined by an ordering $\calo\in\calB_{sg}^{(m)}$ and $\mathbf{U}_k$, $k\geq0$, are UBC transformations, then the following two relations are equivalent:
\begin{itemize}
\item[(i)] $\displaystyle \lim_{k\rightarrow\infty} \off\left(\widehat{\mathbf{H}}_{p(k)q(k)}^{(k+1)}\right)=0,$
\item[(ii)] $\displaystyle \lim_{k\rightarrow\infty} \off\left(\mathbf{H}^{(k)}\right)=0.$
\end{itemize}
\end{Theorem}

\begin{proof}
The proof is similar to the proof of~\cite[Proposition 4.3]{BePe24}. The difference is in the fact that here we have a block process. Therefore, instead of the Jacobi annihilators and operators from~\cite{BeHa21}, one should use their block counterparts from~\cite{BeHa24}. Apart from that, the proof goes the same way. Additionally, the statement of the theorem can be obtained as a special case of~\cite[Theorem 5.1]{Hari15}.
\end{proof}

Furthermore, we prove two propositions.

\begin{Proposition}\label{prop:blockdelta}
Let $\mathbf{A}^{(k)}$, $k\geq 0$, be a sequence generated by the iterative process~\eqref{blockEber_process}. Then, for the reduction of the Frobenius norm in the $k$th step the following inequality holds,
$$\Delta_k=\|\mathbf{A}^{(k)}\|_F^2-\|\mathbf{A}^{(k+1)}\|_F^2 \geq \frac{1}{3} \sum_{l=0}^{L-1} \frac {|\doubletilde{c}_{r_l s_l}^{(l)}|^2}{\|\mathbf{A}\|_F^2}\geq0, \quad \text{for } r_l<s_l,$$
where $C\Big(\doublewidetilde{\mathbf{A}}^{(l)}\Big)=\big(\doubletilde{c}_{ij}^{(l)}\big)$ and $L$ is as defined in~\eqref{indexl}. 
\end{Proposition}

\begin{proof}
After each transformation of the form~\eqref{eq:coreprocess}, according to~\eqref{eq:delta}, we have
$$\triangle_l\geq \frac{1}{3} \frac {|\doubletilde{c}_{r_l s_l}^{(l)}|^2}{\|\mathbf{A}\|_F^2}\geq 0.$$
Therefore, for the reduction of the Frobenius norm in the $k$th step of the process~\eqref{blockEber_process}, regarding the dependence of $l$ on $k$, we have
$$\triangle_k=\|\mathbf{A}^{(k)}\|_F^2-\|\mathbf{A}^{(k+1)}\|_F^2 =\sum_{l=0}^{L-1} \left(\|\doublewidetilde{\mathbf{A}}^{(l)}\|_F^2-\|\doublewidetilde{\mathbf{A}}^{(l+1)}\|_F^2\right) =\sum_{l=0}^{L-1} \triangle_l \geq \frac{1}{3} \sum_{l=0}^{L-1} \frac {|\doubletilde{c}_{r_l s_l}^{(l)}|^2}{\|\mathbf{A}\|_F^2}\geq 0.$$
\end{proof}

\begin{Proposition}\label{prop:Cpq}
Let $\mathbf{A}^{(k)}$, $k\geq 0$, be a sequence generated by the iterative process~\eqref{blockEber_process}. Then,
\begin{equation}\label{eq:limit_Cpq}
\lim_{k\to \infty}{\off\left(\widehat{C}\big(\widetilde{\mathbf{A}}^{(k)}\big)_{pq}\right)}= 0,
\end{equation}
where $\widehat{C}\big(\widetilde{\mathbf{A}}^{(k)}\big)_{pq}$  is the pivot submatrix of ${C}\big(\widetilde{\mathbf{A}}^{(k)}\big)$.
\end{Proposition}

\begin{proof}
From the previous proposition, we see that the sequence $(\|\mathbf{A}^{(k)}\|_F^2$, $k\geq0)$, is non-increasing. Since it is bounded from below, it is convergent.
{That is, $\|\mathbf{A}^{(k)}\|_F^2-\|\mathbf{A}^{(k+1)}\|_F^2\rightarrow0$ and, consequently,  $\|\doublewidetilde{\mathbf{A}}^{(l)}\|_F^2 -\|\doublewidetilde{\mathbf{A}}^{(l+1)}\|_F^2\rightarrow0$. Then, the fact that the left-hand side in~\eqref{eq:delta} tends to zero implies}
\begin{equation}\label{eq:limit_sum_crs}
\lim_{k\to \infty} \sum_{l=0}^{L-1} {|\doubletilde{c}_{r_l s_l}^{(l)}|^2}=0, \quad \text{for } r_l<s_l.
\end{equation}
The notation of the limit in~\eqref{eq:limit_sum_crs} makes sense because the matrices $\doublewidetilde{\mathbf{A}}^{(l)}$, and therefore $C\big(\doublewidetilde{\mathbf{A}}^{(l)}\big)$, $0\leq l<L$, depend on $k$.

It follows from the limit~\eqref{eq:limit_sum_crs} that
$$\lim_{k\to \infty} {\doubletilde{c}^{(l)}_{r_l s_l}}=0, \quad \text {for } 0\leq l<L.$$
The matrix $C\big(\doublewidetilde{\mathbf{A}}^{(l)}\big)$ is Hermitian. Thus, we also have $\lim_{k\to \infty} {\doubletilde{c}^{(l)}_{s_l r_l}}=0$, for $0\leq l<L$.
Now, the assertion~\eqref{eq:limit_Cpq} follows directly from~\eqref{eq:coreprocess} and the fact that the index pairs $(r_l,s_l)$, $0\leq l<L$, correspond to the off-diagonal elements of the pivot submatrix.
\end{proof}

Moreover, we are going to use two results from~\cite{Hari82} for the element-wise Eberlein method~\eqref{Eber-elements}. For $A=A^{(0)}\in\C^{n\times n}$, $k\geq0$, for an arbitrary rotation $R_k$ and $B^{(k)}$ defined as in Section~\ref{sec:prelim}, and $C\big(R_k^*A^{(k)}R_k\big)=C\big(\widetilde{A}^{(k)}\big)=\big(\tilde{c}_{ij}^{(k)}\big)$, the following inequalities hold:
\begin{itemize}
\item[(i)] 
\begin{equation}\label{eq:hariE}
\|A^{(k+1)}-R_k^*A^{(k)}R_k\|_F^2\leq\frac{3}{2}n^2|\tilde{c}_{pq}^{(k)}|,
\end{equation}
\item[(ii)] 
\begin{equation}\label{eq:hariF}
\|B^{(k+1)}-R_k^*B^{(k)}R_k\|_F^2\leq\frac{3}{2}n^2|\tilde{c}_{pq}^{(k)}|.
\end{equation}
\end{itemize}

Now, we are ready to prove the main theorem.

\begin{Theorem}\label{theorem:block_sg}
Let $\mathbf{A}\in\C^{n\times n}$ be a block matrix with partition $\pi=(n_1,\ldots,n_m)$, and let $(\mathbf{A}^{(k)},k\geq0)$  be a sequence generated by the block Eberlein method under a generalized serial pivot strategy defined by an ordering $\calo\in\calB_{sg}^{(m)}$. Let the matrices $\mathbf{B}^{(k)}$ be defined as in~\eqref{eq:blockB}, and the matrices $C\big(\mathbf{A}^{(k)}\big)$ as in~\eqref{eq:C}.
\begin{itemize}
\item[(i)] The sequence of the Hermitian parts of $\mathbf{A}^{(k)}$ converges to a diagonal matrix,
$$\lim_{k\to\infty}\off\big(\mathbf{B}^{(k)}\big)=0.$$
\item[(ii)] The sequence of matrices $\mathbf{A}^{(k)}$ converges to a normal matrix,
$$\lim_{k\to\infty}C\big(\mathbf{A}^{(k)}\big)=0.$$
\item[(iii)] The sequence of the Hermitian parts of $\mathbf{A}^{(k)}$ converges to a fixed diagonal matrix,
$$\lim_{k\to\infty}\mathbf{B}^{(k)}=\diag(\mu_1,\mu_2,\ldots,\mu_n),$$
where $\mu_i$, $1\leq i\leq n$, are real parts of the eigenvalues of $\textbf{A}$.
\item[(iv)] If $\mu_i\neq\mu_j$, then $\lim_{k\to\infty}a_{ij}^{(k)}=0$ and $\lim_{k\to\infty}a_{ji}^{(k)}=0$.
\end{itemize}
\end{Theorem}

\begin{proof}
The proof follows the proof of Theorem 4.4 from~\cite{BePe24}. However, since we work with block matrices, many details must be clarified.

\begin{itemize}
\item[(i)]
Set $\widetilde{\mathbf{B}}^{(k)}=\mathbf{R}_k^*\mathbf{B}^{(k)}\mathbf{R}_k$ and $\mathbf{F}^{(k)}=\mathbf{B}^{(k+1)}-\widetilde{\mathbf{B}}^{(k)}$, $k\geq 0$. 
We will consider the iterative process
\begin{equation}\label{tm:blockB}
\mathbf{B}^{(k+1)}=\mathbf{R}_k^*\mathbf{B}^{(k)}\mathbf{R}_k+\mathbf{F}^{(k)}, \quad k\geq0.
\end{equation}

We observe the impact of the norm-reducing transformation $\textbf{S}_k$ and the iterations~\eqref{eq:coreprocess}. Using the inequality~\eqref{eq:hariF}, we get
\begin{align*}
\|\mathbf{F}^{(k)}\|_F^2 & = \|\mathbf{B}^{(k+1)}-\widetilde{\mathbf{B}}^{(k)}\|_F^2 = \|\doublewidetilde{\mathbf{B}}^{(L)}-\doublewidetilde{\mathbf{B}}^{(0)}\|_F^2 \\
& \leq \sum_{l=0}^{L-1} \|\doublewidetilde{\mathbf{B}}^{(l+1)}
-\doublewidetilde{\mathbf{B}}^{(l)}\|_F^2\leq
\frac{3}{2} n^2 \sum_{l=0}^{L-1} {|\doubletilde{c}_{r_l s_l}^{(l)}|^2}.
\end{align*}
Then, relation~\eqref{eq:limit_sum_crs} from the Proposition~\ref{prop:Cpq} implies
\begin{equation}\label{tm_limF}
\lim_{k\to\infty}\mathbf{F}^{(k)}=0.
\end{equation}
Matrices $\textbf{R}_k$ are assumed to be UBC matrices, so we conclude that the relation~\eqref{tm:blockB} defines a Jacobi-type process like the one in the statement of the Theorem~\ref{tm:blockJacobi}.

Further on, for the pivot submatrices determined by the pivot pair $(p,q)=(p(k),q(k))$, we have
$$\widehat{\mathbf{B}}_{pq}^{(k+1)}=\widehat{\widetilde{\mathbf{B}}}_{pq}^{(k)}+\widehat{\mathbf{F}}_{pq}^{(k)}.$$
It follows from~\eqref{tm_limF} that $\lim_{k\to\infty}\widehat{\mathbf{F}}_{pq}^{(k)}=0$. Since the rotation $\mathbf{R}_k$ is chosen to diagonalize $\widehat{\mathbf{B}}_{pq}^{(k)}$, submatrix $\widehat{\widetilde{\mathbf{B}}}_{pq}^{(k)}$ is diagonal. Therefore, 
$$\lim_{k\to\infty}\off\big(\widehat{\mathbf{B}}_{pq}^{(k+1)}\big)=0.$$
Thus, Theorem~\ref{tm:blockJacobi} implies
\begin{equation*}
\lim_{k\to\infty}{\off\big(\mathbf{B}^{(k)}\big)}=0. 
\end{equation*}

\item[(ii)]
Set $\mathbf{E}^{(k)}=\mathbf{A}^{(k+1)}-\widetilde{\mathbf{A}}^{(k)}$, $k\geq0$. Again, we consider the impact of the transformations $\textbf{S}_k$ as given in~\eqref{eq:coreprocess}. We use the inequality~\eqref{eq:hariE} to obtain
\begin{align*}
\|\mathbf{E}^{(k)}\|_F^2 & = \|\mathbf{A}^{(k+1)}-\widetilde{\mathbf{A}}^{(k)}\|_F^2 = \|\doublewidetilde{\mathbf{A}}^{(L)}-\doublewidetilde{\mathbf{A}}^{(0)}\|_F^2 \\
& \leq \sum_{l=0}^{L-1} \|\doublewidetilde{\mathbf{A}}^{(l+1)}
-\doublewidetilde{\mathbf{A}}^{(l)}\|_F^2\leq
\frac{3}{2} n^2 \sum_{l=0}^{L-1} {|\doubletilde{c}_{r_l s_l}^{(l)}|^2}.
\end{align*}
Hence, relation~\eqref{eq:limit_sum_crs} implies
\begin{equation}\label{tm_limE}
\lim_{k\to\infty}\mathbf{E}^{(k)}=0.
\end{equation}

Now, regarding the distance from normality, we have
$$C\big(\mathbf{A}^{(k+1)}\big)=C\big(\widetilde{\mathbf{A}}^{(k)}+\mathbf{E}^{(k)}\big), \quad k\geq0.$$
From the definition of the operator $C$, after a short calculation we obtain
\begin{equation}\label{tm:CAblock}
C\big(\mathbf{A}^{(k+1)}\big) = C\big(\widetilde{\mathbf{A}}^{(k)}\big)+\mathbf{W}^{(k)}, \quad k\geq0
\end{equation}
where
$$\mathbf{W}^{(k)}=\mathbf{A}^{(k+1)}\big(\mathbf{E}^{(k)}\big)^*-\big(\mathbf{A}^{(k+1)}\big)^*\mathbf{E}^{(k)}+\mathbf{E}^{(k)}\big(\widetilde{\mathbf{A}}^{(k)}\big)^*-\big(\mathbf{E}^{(k)}\big)^*\widetilde{\mathbf{A}}^{(k)},$$
and
$$\|\mathbf{W}^{(k)}\|_F \leq 2\|\mathbf{E}^{(k)}\|_F\left(\|\mathbf{A}^{(k+1)}\|_F+\|\widetilde{\mathbf{A}}^{(k)}\|_F\right).$$
Proposition~\ref{prop:blockdelta} implies
$$\|\mathbf{W}^{(k)}\|_F \leq 4\|\mathbf{E}^{(k)}\|_F\|\widetilde{\mathbf{A}}^{(k)}\|_F \leq 4\|\mathbf{E}^{(k)}\|_F\|\mathbf{A}\|_F.$$
Then, it follows from the relation~\eqref{tm_limE} that
\begin{equation}\label{tm:blockW}
\lim_{k\to\infty}\|\mathbf{W}^{(k)}\|_F=0.
\end{equation}

For $C\big(\widetilde{\mathbf{A}}^{(k)}\big)$, from the definition of the operator $C$, it follows
$$C\big(\widetilde{\mathbf{A}}^{(k)}\big)=C\big(\mathbf{R}_k^*\mathbf{A}^{(k)}\mathbf{R}_k\big)=\mathbf{R}_k^*C\big(\mathbf{A}^{(k)}\big)\mathbf{R}_k, \quad k\geq0.$$
This means that the iterations~\eqref{tm:CAblock} can be written as
\begin{equation}\label{tm:Cprocess}
C\big(\mathbf{A}^{(k+1)}\big)=\mathbf{R}_k^*C\big(\mathbf{A}^{(k)}\big)\mathbf{R}_k+\mathbf{W}^{(k)}, \quad k\geq0.
\end{equation}
It is easy to check that matrix $C\big(\mathbf{A}^{(0)}\big)=C\big(\mathbf{A}\big)$ is Hermitian. Matrices $\mathbf{R}_k$ are UBC transformations and~\eqref{tm:blockW} holds. Therefore, relation~\eqref{tm:Cprocess} defines a Jacobi-type process satisfying the statement of the Theorem~\ref{tm:blockJacobi}.
For the pivot submatrices in the $k$th step we have
$$\widehat{C}\big(\mathbf{A}^{(k+1)}\big)_{pq}=\widehat{C}\big(\widetilde{\mathbf{A}}^{(k)}\big)_{pq}+\widehat{\mathbf{W}}_{pq}^{(k)}.$$
Relations~\eqref{eq:limit_Cpq} and~\eqref{tm:blockW} imply $\lim_{k\to\infty}\off\big(\widehat{C}(\mathbf{A}^{(k+1)})_{pq}\big)=0$. Thus,
\begin{equation}\label{tm:offC_block}
\lim_{k\to\infty}\off\big(C\big(\mathbf{A}^{(k)}\big)\big)=0,
\end{equation}
follows from the Theorem~\ref{tm:blockJacobi}.

So far, we showed that the off-diagonal elements of $C\big(\mathbf{A}^{(k)}\big)$ converge to zero. We should show that the diagonal elements have the same property.
We can write $C\big(\mathbf{A}^{(k)}\big)$, $k\geq0$, as a sum of its Hermitian part $\mathbf{B}^{(k)}$ and skew-Hermitian part $\mathbf{Z}^{(k)}$. Then,
\begin{align}
C\big(\mathbf{A}^{(k)}\big) & = \big(\mathbf{B}^{(k)}+\mathbf{Z}^{(k)}\big)\big(\mathbf{B}^{(k)}+\mathbf{Z}^{(k)}\big)^*-\big(\mathbf{B}^{(k)}+\mathbf{Z}^{(k)}\big)^*\big(\mathbf{B}^{(k)}+\mathbf{Z}^{(k)}\big) \nonumber\\
& = 2\big(\mathbf{Z}^{(k)}\mathbf{B}^{(k)}-\mathbf{B}^{(k)}\mathbf{Z}^{(k)}\big). \label{tm:CBG_block}
\end{align}
The diagonal elements of $C\big(\mathbf{A}^{(k)}\big)$ are then given by
$$c_{ii}^{(k)}=2\sum_{j=1}^n \left(z_{ij}^{(k)} b_{ji}^{(k)}-b_{ij}^{(k)}z_{ji}^{(k)}\right), \quad 1\leq i\leq n.$$
Part $(i)$ of this theorem implies 
$$\lim_{k\to \infty} b_{ij}^{(k)}=0, \quad \text{for } i\neq j,$$
that is,
$$\lim_{k\to\infty}c_{ii}^{(k)}=2\left(z_{ii}^{(k)}b_{ii}^{(k)}-b_{ii}^{(k)}z_{ii}^{(k)}\right)=0.$$
Together with the relation~\eqref{tm:offC_block}, this proves part (ii).

\item[(iii)]
The assertion can be shown the same way as in~\cite[Theorem 4.3]{PupHari98}, using the assertion (ii) of this theorem.

\item[(iv)]
The assertion can be shown the same way as for the element-wise method in~\cite[Theorem 4.4]{BePe24}, using the relation~\eqref{tm:CBG_block} and parts (i)--(iii) of this theorem.
\end{itemize}
\end{proof}

Let us recapitulate the results of Theorem~\ref{theorem:block_sg}. 
\begin{itemize}
\item Starting with an ${n\times n}$ complex block matrix $\textbf{A}^{(0)}=\textbf{A}$ with partition $\pi$ as given in~\eqref{partition_pi}, the sequence of block matrices $(\textbf{A}^{(k)},k\geq 0)$ generated by the block Eberlein method~\eqref{blockEber_process}, under any generalized serial block pivot strategy, converges to a normal matrix $\Lambda$. 
\item If all the eigenvalues of $\mathbf{A}$ have different real parts, then $\Lambda$ is a diagonal matrix with the eigenvalues of $\mathbf{A}$ on the diagonal of $\Lambda$.
\item If there are eigenvalues of $\mathbf{A}$ with the same real parts, then $\Lambda$ is permutation-similar to a block diagonal matrix with the block sizes equal to the number of times the same real part appears in the spectrum of $\textbf{A}$. 
These diagonal blocks do not necessarily match the partition $\pi$.
\item In the case where there are repeating real parts in the spectrum of $\mathbf{A}$, the eigenvalues with non-repeating real parts can be read from the diagonal of $\Lambda$, while the other eigenvalues can be obtained from the blocks of $\Lambda$, which comes down to solving the eigenvalue problems of the (small) matrix blocks.
\item The Hermitian part of $\textbf{A}^{(k)}$ always converges to a diagonal matrix with the real parts of the eigenvalues of $\textbf{A}$ on the diagonal. 
\end{itemize}

In our numerical tests, we have observed some interesting things related to the blocks of $\Lambda$.
Blocks do not appear, in practice, in the case of multiple complex eigenvalues with the same real and the same imaginary part. This includes the cases of multiple purely real or purely imaginary eigenvalues. In practice, blocks appear only if there are complex eigenvalues with the same real, but different imaginary parts.
According to this observation, it is useful to precondition the starting matrix $\textbf{A}$.

It is easy to check that, for a random $0\neq d\in\C$, if $(\lambda,x)$ is an eigenpair of a matrix $M$, then $(d\lambda,x)$ is an eigenpair of $dM$. We take $d\in\C$ such that $\text{Im}(d)\neq0$ and apply the block Eberlein method to $d\textbf{A}$. Then, with probability one, matrix $d\textbf{A}$ does not have eigenvalues with the same real and different imaginary parts. Therefore, the block Eberlein algorithm applied on $d\textbf{A}$, in practice, results in a diagonal matrix $\Lambda_d$ with the eigenvalues of $d\textbf{A}$ on the diagonal. Then, the eigenvalues of $\textbf{A}$ are obtained by dividing the diagonal elements of $\Lambda_d $ by $d$.

\section{Numerical examples}\label{sec:numerical}

{Before we go through numerical examples, let us look at the complexity of the element-wise and block variants of the Eberlein algorithm. One iteration of the element-wise Eberlein method consists of two parts. The first is finding $R_k$ and transforming $A^{(k)}\to \widetilde{A}^{(k)}$, while the second part is finding $S_k$ and updating $\widetilde{A}^{(k)}\to A^{(k+1)}$. Number of operations for computing those is $c_1n+c_2$ and $c_3n+c_4$, respectively, for some integers $c_i$, $i=1,2,3,4$. That is, complexity of a single step is $O(n)$. Overall complexity of one sweep of the element-wise method is
$$O(n)\cdot\frac{n(n-1)}{2}=O(n^3).$$
One iteration of our block Eberlein method also consists of two parts, described by~\eqref{eq:blockpart1} and~\eqref{eq:blockpart2}. The first part  diagonalizes an $(n_p+n_q)\times(n_p+n_q)$ matrix which comes with complexity $O((n_p+n_q)^3)$. The operation count for the second part is $(c_1n+c_2)\cdot(n_p+n_q)(n_p+n_q-1)/2,$ for some integers $c_1,c_2$. Hence, complexity of the second part is $O(n(n_p+n_q)^2)$. Then, assuming that $n_p+n_q\ll n$, complexity of one iteration corresponds to $O(n)$. Complexity of one sweep of the block method is
$$O(n)\cdot\frac{m(m-1)}{2}=O(nm^2),$$
where $m<n$ is the number of blocks.}

Now we present numerical tests for the Algorithm~\ref{agm:blockeberlein}. We use the row-wise pivot strategy which belongs to the set of the generalized serial pivot strategies. All experiments were performed in Matlab R2024b.

In accordance with Theorem~\ref{theorem:block_sg}, we are interested in the behavior of the off-norms of the matrices $\textbf{A}^{(k)}$ and $\textbf{B}^{(k)}$, and the norm of the matrices $C(\textbf{A}^{(k)})$, $k\geq 0$. The results are presented in logarithmic scale. We stop the algorithm when the change in the off-norm of $\textbf{B}^{(k)}$ is smaller than the tolerance, in our case $10^{-10}$. 
We take the partition $\pi=(n_1,n_2,\ldots,n_m)$ to have all blocks of the same size, $n_1=n_2=\ldots =n_m$. We tested the algorithm for different block sizes, five, ten, and twenty. Each line in the figures represents the results for a different block size. 

Additionally, we test the accuracy of the block Eberlein method. We show the relative errors in the real and imaginary parts of the diagonal elements of $\Lambda$, with respect to the eigenvalues obtained by the Matlab function $\verb|eig|$. Moreover, the columns $t_i$, $1\leq i\leq n$, of the matrix $\textbf{T}_K$ acquired by the Algorithm~\ref{agm:blockeberlein} represent the eigenvectors corresponding to the eigenvalue $\lambda_i$. To depict the accuracy of the computed eigenvectors, we look at the values of 
\begin{equation}\label{num_eigvec}
|\mathbf{A}t_i-\lambda_it_i|, \quad 1\leq i\leq n.
\end{equation}

{We start our numerical tests with a diagonalizable matrix. For $n=500$, we construct a test matrix $\textbf{A}_0\in \C^{n\times n}$ as follows:}
{
\begin{itemize}
\item $\verb|d=randn(n,1)+1i*randn(n,1); |$
\item $\Sigma=\verb|diag(d);|$
\item $\verb|Q=orth(randn(n)+1i*randn(n));|$
\item $\verb|A_0=Q*|\Sigma\verb|*Q';|$
\end{itemize}}

\begin{figure}[ht]
\begin{subfigure}[t]{\textwidth}
\centering
\includegraphics[width=0.61\textwidth]{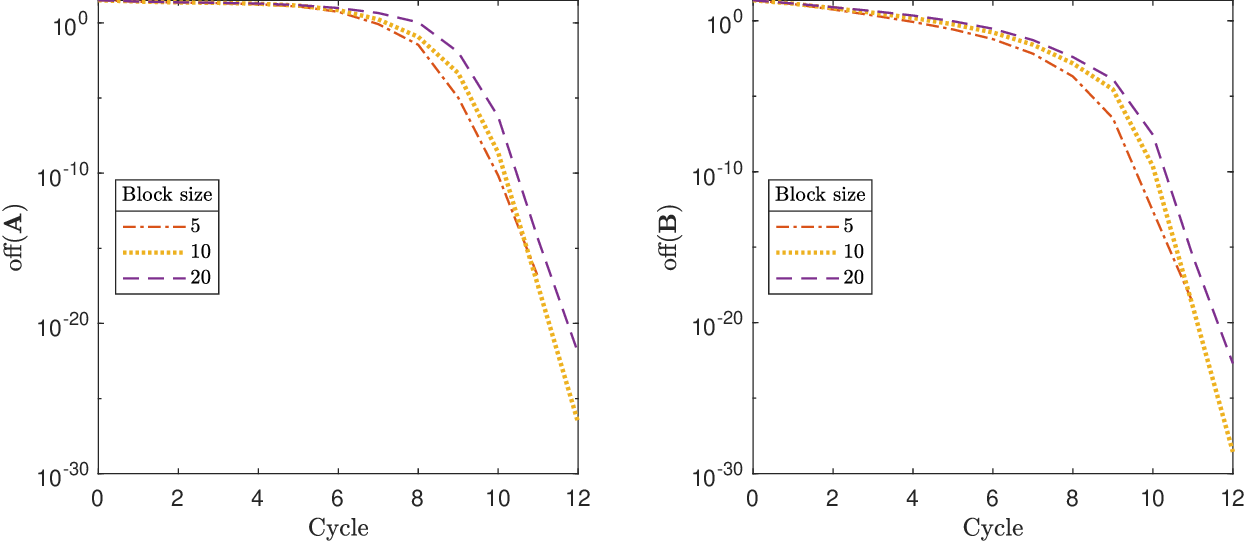}
\caption{Change in $\off(\textbf{A}^{(k)})$ and $\off(\textbf{B}^{(k)})$ for different block sizes.}\label{fig:A0_cyc}
\end{subfigure}
\begin{subfigure}[t]{\textwidth}
\centering
\includegraphics[width=0.9\textwidth]{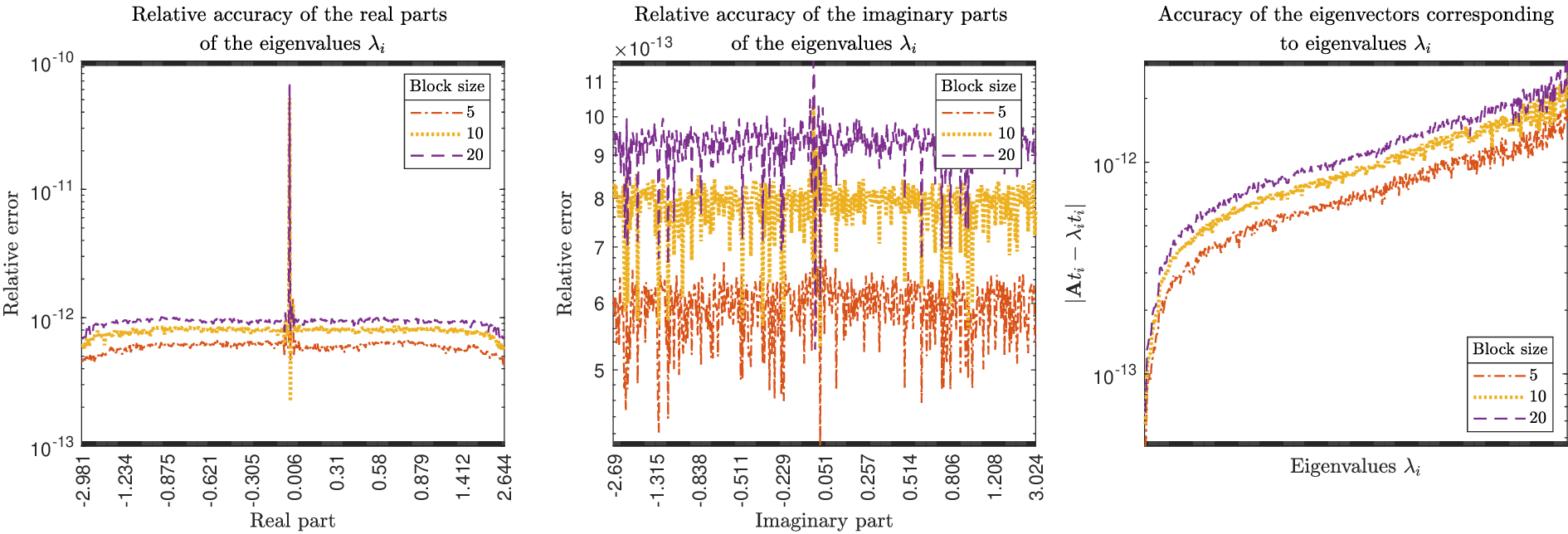}
\caption{Accuracy of the eigenvalues in comparison to the Matlab eig function, and accuracy of the eigenvectors.}\label{fig:A0_acc}
\end{subfigure}
\caption{Results for the test matrix $\textbf{A}_0$ with $n=500$.}\label{fig:A0}
\end{figure}

{We apply the block Eberlein algorithm to $\textbf{A}_0$ for different block sizes and display the outcome in Figure~\ref{fig:A0}. 
We do not observe the distance to normality,  $\|C(\textbf{A}^{(k)})\|_F$, since $\textbf{A}_0$ is diagonalizable and normal. This essentially means that all transformations $\textbf{S}_k$ are equal to identity. After a sufficient but not very large number of cycles, both $\textbf{A}^{(k)}$ and its Hermitian part are diagonalized, as seen in Figure~\ref{fig:A0_cyc}. Accuracy of the real and imaginary parts of the obtained eigenvalues, as well as the accuracy of the corresponding eigenvectors are shown in Figure~\ref{fig:A0_acc}.}

\begin{figure}[ht]
\begin{subfigure}[t]{\textwidth}
\centering
\includegraphics[width=0.9\textwidth]{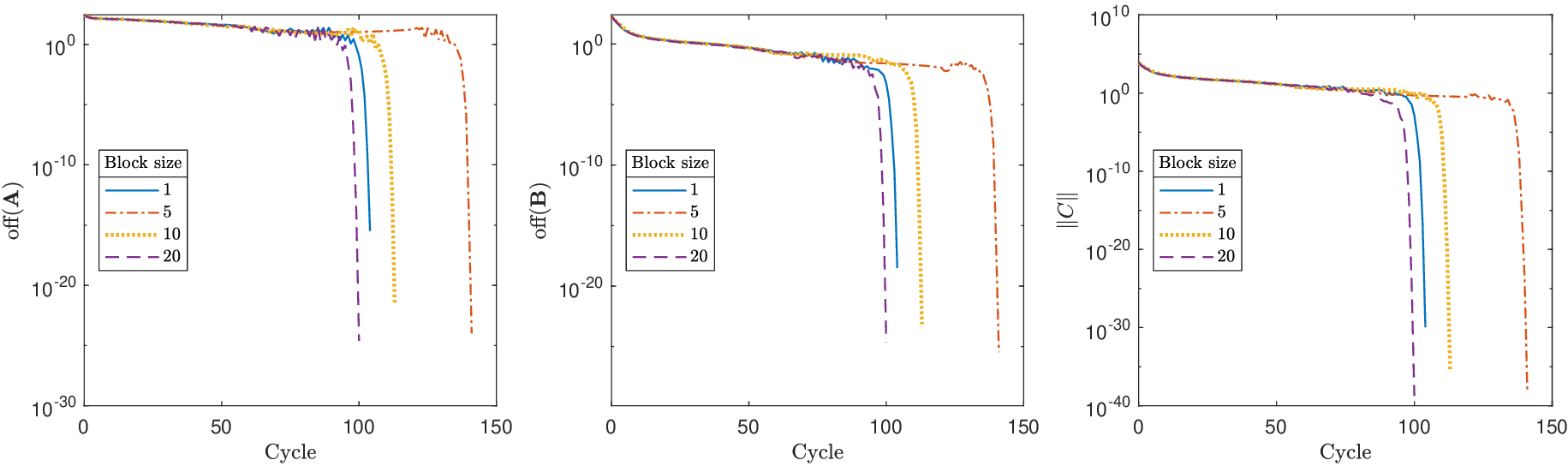}
\caption{Change in $\off(\textbf{A}^{(k)})$, $\off(\textbf{B}^{(k)})$, and $\|C(\textbf{A}^{(k)})\|_F$ for different block sizes.}\label{fig:A1_cyc}
\end{subfigure}
\begin{subfigure}[t]{\textwidth}
\centering
\includegraphics[width=0.61\textwidth]{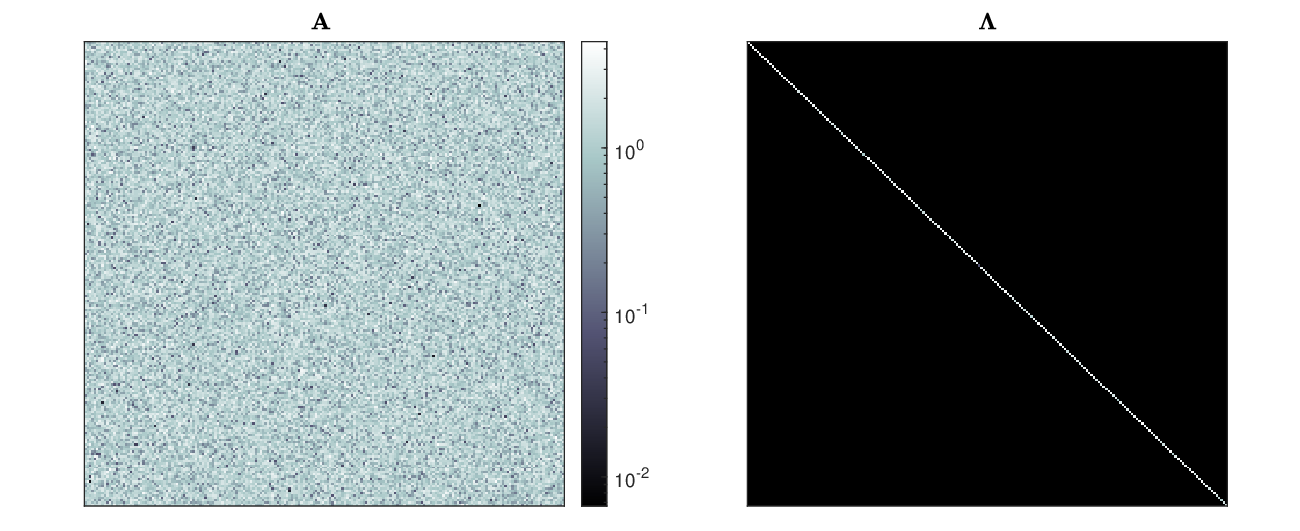}
\caption{Matrix structure.}\label{fig:A1_heat}
\end{subfigure}
\bigskip
\begin{subfigure}[t]{\textwidth}
\centering
\includegraphics[width=0.9\textwidth]{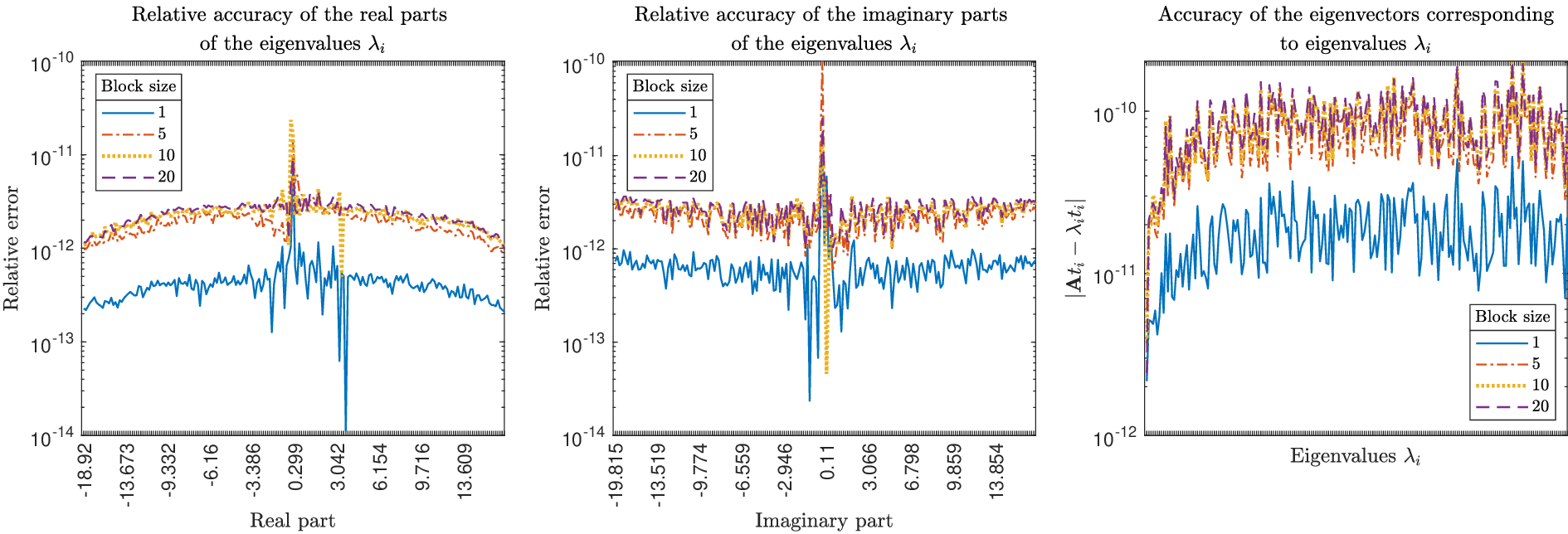}
\caption{Accuracy of the eigenvalues in comparison to the Matlab eig function, and accuracy of the eigenvectors.}\label{fig:A1_acc}
\end{subfigure}  
\caption{Results for the test matrix $\textbf{A}_1$ with $n=200$.}\label{fig:A1}
\end{figure}

Moreover, we demonstrate the algorithm on a random matrix. For $n=200$, we construct the test matrix $\textbf{A}_1\in\C^{n\times n}$ as:
\begin{itemize}
\item  $\verb|A_1=randn(n)+1i*randn(n);|$
\end{itemize}  
Generically, random matrices are not normal and have different eigenvalues. Figure~\ref{fig:A1} shows the results of the block Eberlein method applied to $\textbf{A}_1$. As expected, in Figure~\ref{fig:A1_cyc} we see that $\off(\textbf{B}^{(k)})$ and $\|C(\textbf{A}^{(k)})\|_F$, $k\geq 0$, converge to zero. Since the eigenvalues of $\textbf{A}_1$ are simple, $\off(\textbf{A}^{(k)})$, $k\geq0$, converges to zero, as well. 
{Compared to the previous example, even though the matrix is smaller, it took significantly more cycles until the algorithm converged. This is because the matrix $\textbf{A}_1$ is not normal. Related to the discussion of complexity of the block algorithm from the beginning of this section, we can see that the efficiency of the method depends a lot on well-chosen block sizes.}

In Figure~\ref{fig:A1_heat} we have the structure of the starting matrix $\textbf{A}_1$ and of the resulting matrix $\Lambda$ obtained by the Algorithm~\ref{agm:blockeberlein}, for a block size equal to $20$. In particular, we observe the logarithm of the absolute values of the elements of $\textbf{A}_1$ and $\Lambda$. The elements that are smaller in absolute value are in darker shades. The obtained matrix $\Lambda$ is diagonal. 

In Figure~\ref{fig:A1_acc} we see the relative accuracy of the block Eberlein method compared to the Matlab function $\verb|eig|$. Since $\Lambda$ is diagonal, it carries the approximations of the eigenvalues of $\mathbf{A}_1$. On the $x$-axes we have real (imaginary) parts of the eigenvalues obtained by $\verb|eig|$, arranged in increasing order. The relative errors, in both real and imaginary parts, of the obtained eigenvalues are close to $10^{-12}$, for all block sizes. 
The third graph in this figure shows the accuracy of the eigenvectors, as given in~\eqref{num_eigvec}, for different block sizes. On the $x$-axis, the computed eigenvalues are arranged in ascending order, with respect to the absolute value. In general, using the partition with smaller sized blocks yields more accurate approximations.

Next, we test our algorithm on a matrix with the eigenvalues that have the same real parts. The spectrum of the test matrix $\textbf{A}_2\in\C^{n\times n}$ consists of a random complex number $a_1$ and four pairs of complex conjugate numbers $a_i$ and $a_i^*$, with different multiplicities $m_i$, $i=2,3,4,5$. The multiplicities of the eigenvalues add up to $n$, that is, $m_1+2m_2+2m_3+2m_4+2m_5=n$. We consider $\textbf{A}_2$ for $n=200$, $m_1=40$, $m_2=m_3=m_4=m_5=20$. It is constructed as follows:
\begin{itemize}
\item $\verb|a_1=randn(1)+1i*randn(1); a_2=randn(1,4)+1i*randn(1,4);|$
\item $\verb|a=[a_1, a_2, conj(a_2)];|$
\item $\verb|m=[m_1, m_2, m_3, m_4, m_5, m_2, m_3, m_4, m_5];|$
\item $\verb|a=repelem(a, m);|$
\item $\Sigma=\verb|diag(a);|$
\item $\verb|Q=orth(randn(n)+1i*randn(n));|$
\item $\verb|A_2=Q*|\Sigma\verb|*Q';|$
\end{itemize}

\begin{figure}[ht]
\begin{subfigure}[t]{\textwidth}
\centering
\includegraphics[width=0.61\textwidth]{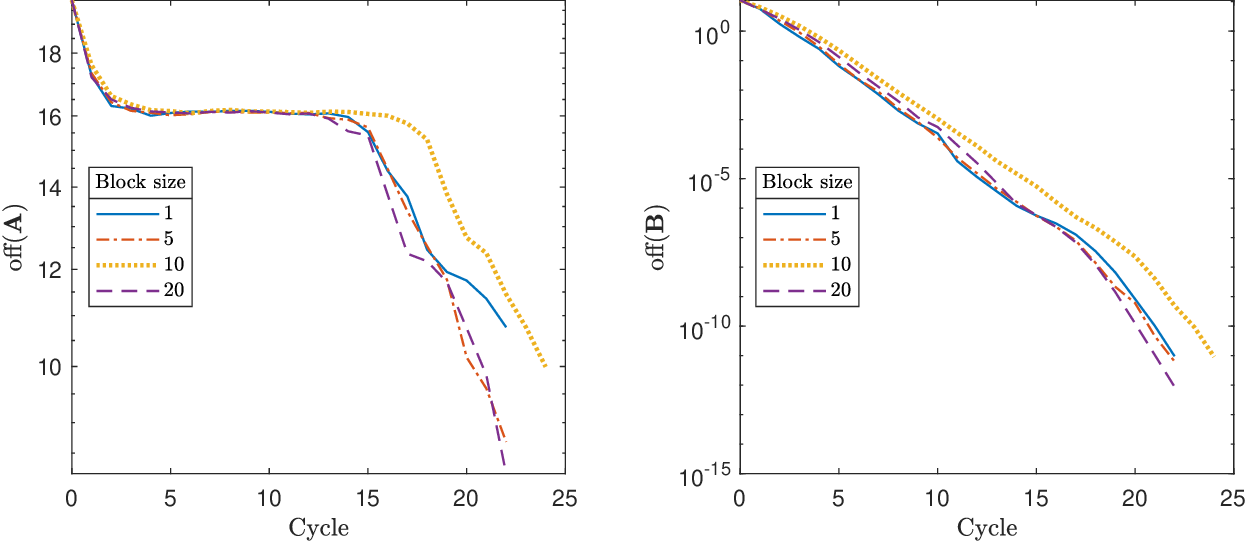}
\caption{Change in $\off(\textbf{A}^{(k)})$ and $\off(\textbf{B}^{(k)})$ for different block sizes.}\label{fig:A2_cyc}
\end{subfigure}
\begin{subfigure}[t]{\textwidth}
\centering
\hspace{0.5cm}\includegraphics[width=0.66\textwidth]{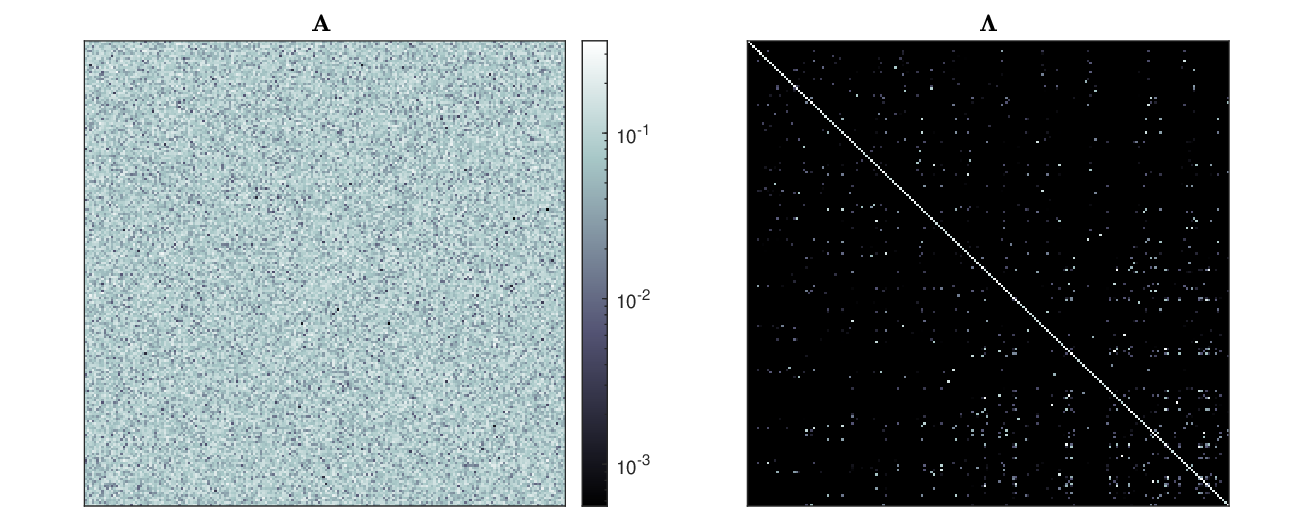}
\caption{Matrix structure, block size = 20.}\label{fig:A2_heat}
\end{subfigure}
\caption{Results for the test matrix $\textbf{A}_2$ with $n=200$, $m_1=40$, $m_2=m_3=m_4=m_5=20$.}\label{fig:A2}
\end{figure}

The results are presented in Figure~\ref{fig:A2}. We can see from Figure~\ref{fig:A2_cyc} that $\off(\textbf{A}^{(k)})$, $k\geq 0$, does not converge to zero. Since some of the eigenvalues of $\mathbf{A}_2$ have the same real part, this corresponds to the results of Theorem~\ref{theorem:block_sg}. The nonzero off-diagonal elements in the resulting matrix $\Lambda$, shown in Figure~\ref{fig:A2_heat}, correspond to the pairs of the complex conjugate eigenvalues $a_i$ and $a_i^*$, $i=2,3,4,5$. The repeating eigenvalue $a_1$ did not create a block and it appears on the diagonal. Note that, for $a_1$, there are no other eigenvalues with the same real but different imaginary part. Despite the repeating eigenvalues, we can see from the Figure~\ref{fig:A2_cyc} that
$\off(\textbf{B}^{(k)})$, $k\geq 0$, converges to zero. Matrix $\textbf{A}_2$ is normal by construction and it stays normal, so there is no need to observe $C(\textbf{A}^{(k)})$. 

\begin{figure}[ht]
\begin{subfigure}[t]{\textwidth}
\centering
\includegraphics[width=0.61\textwidth]{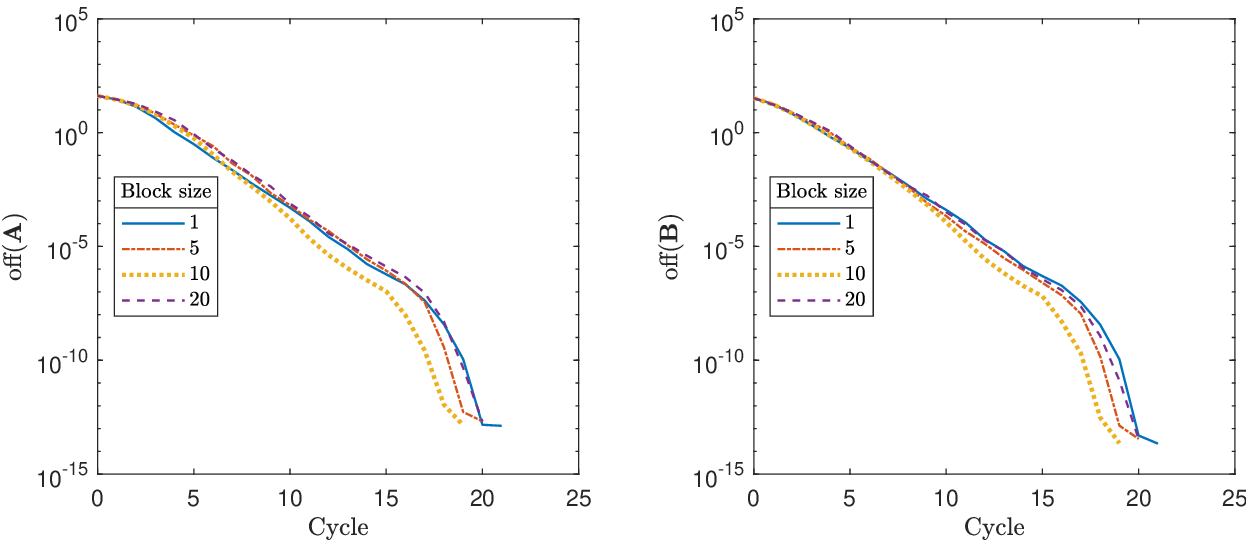}
\caption{Change in $\off(d\textbf{A}^{(k)})$ and $\off(d\textbf{B}^{(k)})$ for different block sizes.}\label{fig:dA2_cyc}
\end{subfigure}
\begin{subfigure}[t]{\textwidth}
\centering
\hspace{0.5cm}\includegraphics[width=0.66\textwidth]{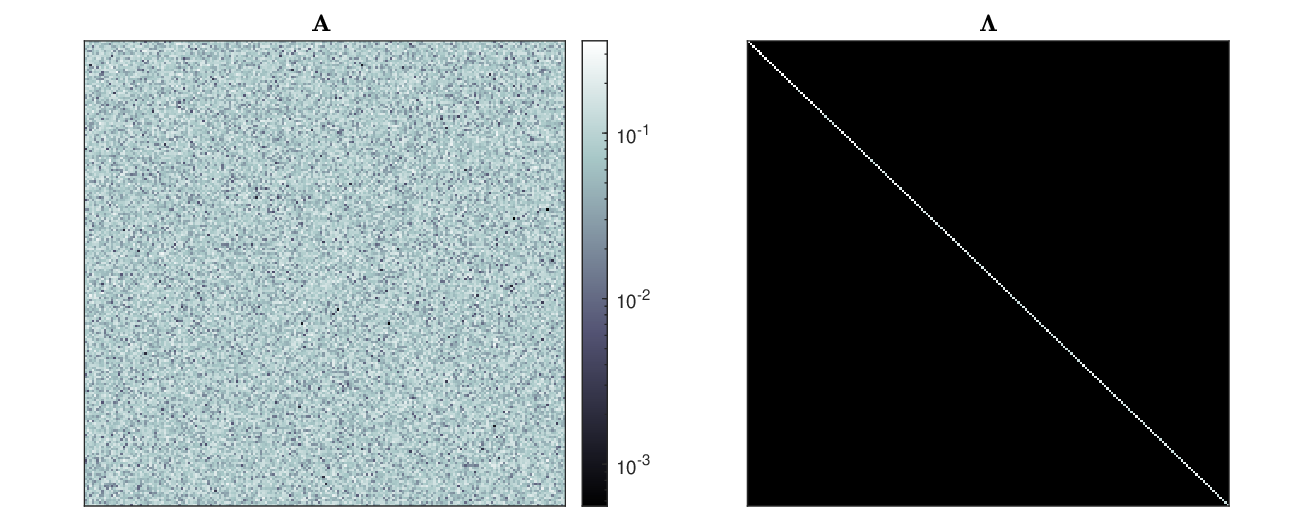}
\caption{Matrix structure, block size = 20.}\label{fig:dA2_heat}
\end{subfigure}
\bigskip
\begin{subfigure}[t]{\textwidth}
\centering
\includegraphics[width=0.9\textwidth]{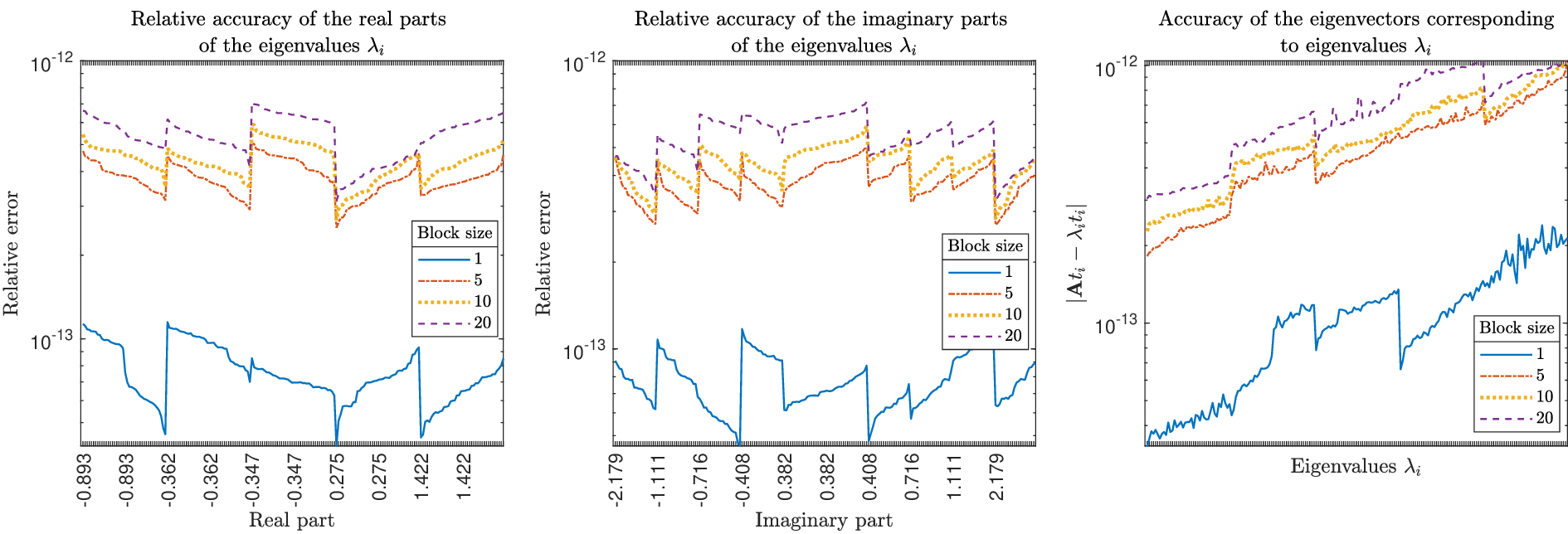}
\caption{Accuracy of the eigenvalues in comparison to the Matlab eig function, and accuracy of the eigenvectors.}\label{fig:dA2_acc}
\end{subfigure}
\caption{Results for the test matrix $\textbf{A}_2$ with $n=200$, $m_1=40$, $m_2=m_3=m_4=m_5=20$, with preconditioning.}\label{fig:dA2}
\end{figure}

We solve this issue and avoid the discussion about the repeating real parts of the eigenvalues by preconditioning the starting matrix. 
We multiply the starting matrix $\textbf{A}_2$ by a complex number $d$, such that Im$(d)\neq0$. Preconditioning step:
\begin{itemize}
\item $\verb|d=randn(1)+1i*randn(1);|$
\item $\verb|dA_2=d*A_2|;$
\end{itemize}
Applying the block Eberlein method to $d\textbf{A}_2$ yields a fully diagonal matrix $\Lambda_d$, as seen in the Figures~\ref{fig:dA2_cyc} and~\ref{fig:dA2_heat}. The eigenvalues of $\textbf{A}_2$ are retrieved by dividing the values on the diagonal of $\Lambda_d$ by $d$. According to Figure~\ref{fig:dA2_acc}, both real and imaginary parts of all the eigenvalues are highly accurate, with respect to the eigenvalues of $\textbf{A}$ acquired by the Matlab function $\verb|eig|$. 
Moreover, Figure~\ref{fig:dA2_acc} shows the accuracy of the eigenvectors obtained by the block Eberlein method with preconditioning. The described simple preconditioning can be  performed before the Eberlein algorithm is applied, as shown in our example. Alternatively, after the Eberlein algorithm has converged to a non-diagonal matrix, one can multiply the resulting matrix with $d$, Im$(d)\neq0$, and run the Eberlein algorithm on a newly obtained matrix.

\begin{figure}[ht]
\begin{subfigure}[t]{\textwidth}
\centering
\includegraphics[width=0.9\textwidth]{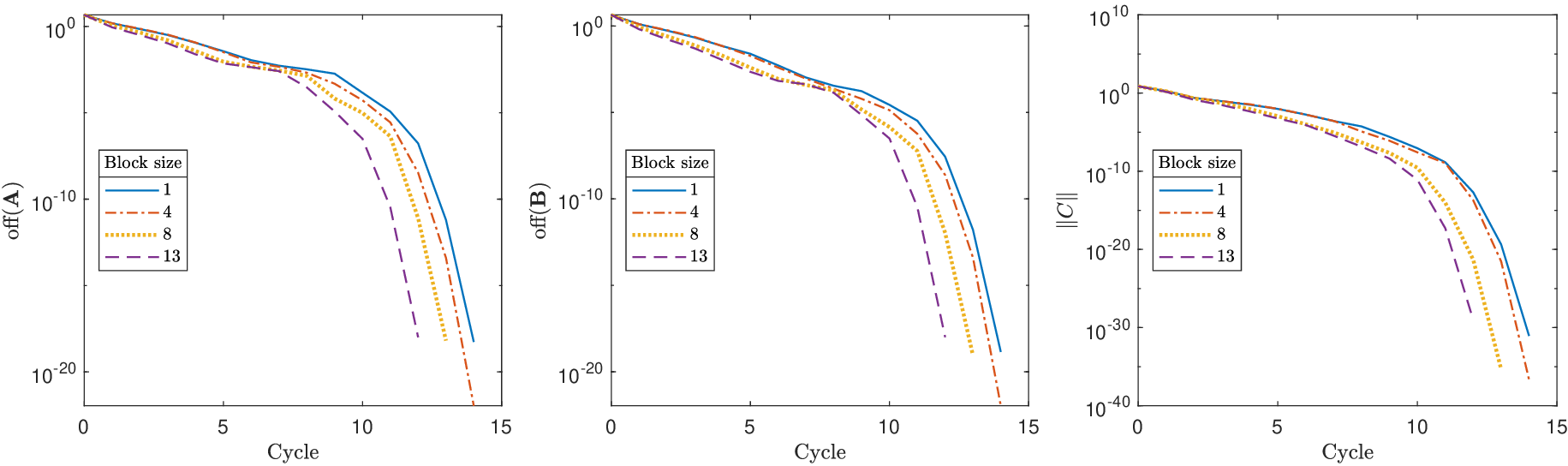}
\caption{Change in $\off(d\textbf{A}^{(k)})$ and $\off(d\textbf{B}^{(k)})$ for different block sizes.}\label{fig:dA3_cyc}
\end{subfigure}
\begin{subfigure}[t]{\textwidth}
\centering
\hspace{0.5cm}\includegraphics[width=0.66\textwidth]{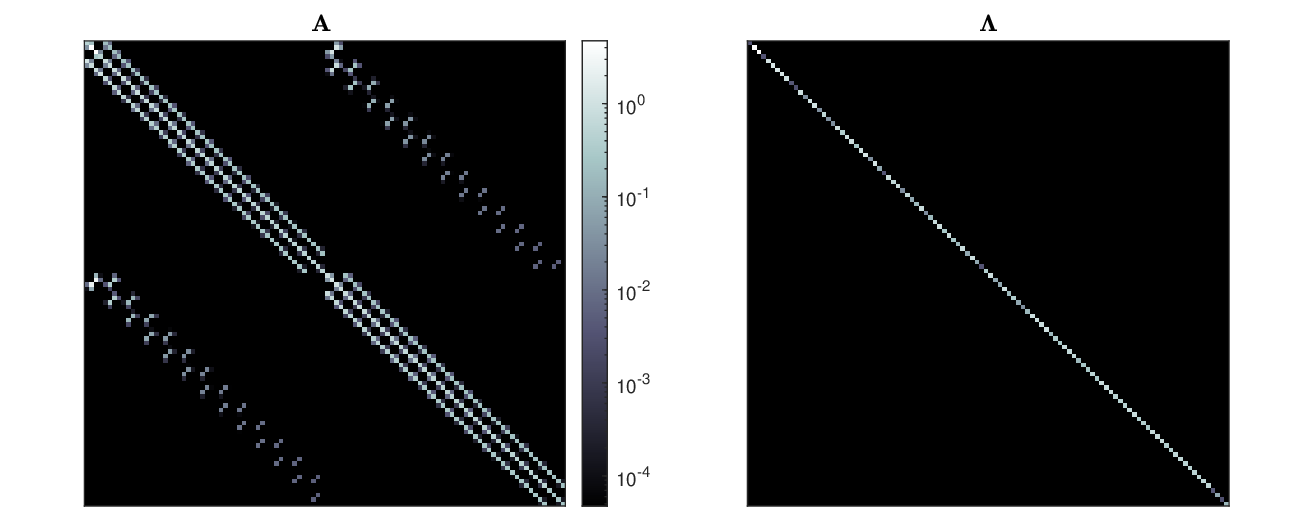}
\caption{Matrix structure, block size = 4.}\label{fig:dA3_heat}
\end{subfigure}
\caption{Results for the test matrix $\textbf{A}_3$, with preconditioning.}\label{fig:dA3}
\end{figure}

{
In the last numerical example, we consider test matrix $\textbf{A}_3=\textup{CK}104\in\R^{104\times 104}$ from~\cite{Bai1997, Matrixmarket}. The majority of the eigenvalues are real, but block Eberlein method reveals several repeating complex conjugate pairs. To obtain all eigenvalues we do a preconditioning step $dA_3$, for some random complex number $d$. The results are given in Figure~\ref{fig:dA3}.}

\section{Conclusion}\label{sec:conclusion}

We presented a block version of the Eberlein diagonalization method. This is, to the best of our knowledge, the first block variant of the Eberlein algorithm. We proved the global convergence of our block algorithm. The convergence results are in line with those for the element-wise method. If all eigenvalues of the starting block matrix $\mathbf{A}$ have different real parts, then the sequence of the matrices obtained by the block Eberlein method converges to a diagonal matrix. Otherwise, it converges to a matrix which is permutation-similar to a block diagonal matrix, with block sizes equal to the number of times the same real part appears in the spectrum of $\mathbf{A}$. In practice, the case of the repeating real parts can be simply solved by preconditioning.

\section*{Acknowledgments}
{The authors thank the anonymous referees for their useful suggestions that improved the paper.}

\bibliographystyle{siam}
\bibliography{blockEberlein.bib}

\end{document}